\documentclass[a4paper,12pt]{article}
\usepackage{geometry}
\usepackage{amssymb}
\usepackage{xypic}
\usepackage[all,2cell]{xy}
\usepackage{graphics}
\usepackage{amsmath}
\usepackage{fullpage}
\usepackage{theorem}
\usepackage{amstext}
\usepackage{mathrsfs}
\usepackage{calligra}
\usepackage{xcolor, soul}
\sethlcolor{yellow}

\newtheorem{theorem}{Theorem}[section]

\theorembodyfont{\upshape}

\newenvironment{proof}[1][Proof]{\textbf{#1.} }{\ \rule{0.5em}{0.5em}}
\begin{document}
	
	\title{On a connection between fuzzy subgroups and $F$-inverse covers of inverse monoids}
	\author{Elton Pasku \\
		Universiteti i Tiran\"es \\
		Fakulteti i Shkencave Natyrore \\
		Departamenti i Matematik\"es\\ 
		Tiran\"e, Albania \\
		elton.pasku@unitir.edu.al}
	
	\date{}
	
	\maketitle

\begin{abstract}
We define two categories, the category $\mathfrak{F}\mathfrak{G}$ of fuzzy subgroups, and the category $\mathfrak{F}\mathfrak{C}$ of $F$-inverse covers of inverse monoids, and prove that $\mathfrak{F}\mathfrak{G}$ fully embeds into $\mathfrak{F}\mathfrak{C}$. This shows that, at least from a categorical viewpoint, fuzzy subgroups belong to the standard mathematics as much as they do to the fuzzy one. \newline
AMS Mathematics  Subject Classification $(2010)$: 03E72, 20N25, 20M18. 
\end{abstract}

\section{Introduction and preliminaries}

The theory of fuzzy sets originates with the article \cite{Z} of Zadeh and has aimed since than to help other branches of mathematics that study ambiguity or uncertainty. Along with fuzzy sets, fuzzy analogues have been developed, in particular the theory of fuzzy groups which started with the paper \cite{R} of Rosenfeld. Given a set $X$, a fuzzy subset $A$ of $X$ is a function $A: X \rightarrow [0,1]$. For every $x \in X$, the value $A(x)$ represents the degree of membership of $x$ in $A$. This is what makes $A$ to look like an uncertain set. On the other hand, the definition of fuzzy groups is a bit more complex and is given below. Let $G$ be a group. A fuzzy subgroup of $G$ is a map $\mu: G \rightarrow [0,1]$ such that:
\begin{enumerate}
	\item [(i)] for all $x,y \in G$, $\mu(x y) \geq \text{min}\{\mu(x),\mu(y)\}$, and
	\item[(ii)] for all $x \in G$, $\mu(x^{-1}) \geq \mu(x)$.
\end{enumerate}	
It turns out that for all $x \in G$, $\mu(x)=\mu(x^{-1})$, and $\mu(e) \geq \mu(x)$ where $e$ is the unit of $G$. There is no restriction if we replace $[0,1]$ in this definition by $\mu([0,1])$, and in this way the definition may be restated as follows. A fuzzy subgroup of $G$ is a triple $(G, \mu, U)$ where $G$ is a group, $U \subseteq [0,1]$ and $\mu: G \rightarrow U$ is a surjective map satisfying the properties:
\begin{enumerate}
	\item [(1)] for all $x,y \in G$, $\mu(x y) \geq \text{min}\{\mu(x),\mu(y)\}$,  and
	\item[(2)] for all $x \in G$, $\mu(x^{-1}) = \mu(x)$.
\end{enumerate}
There is a similarity between the definition of a fuzzy subgroup as a triple, with the definition of a dual premorphism from a group to an inverse monoid as defined in \cite{ASZ}, apart from the fact that here $U$ is not given an inverse monoid structure. Luckily, we can overcome this difficulty very easily. The fact that for all $x \in G$, $\mu(x) \leq \mu(e)$ says exactly that $\text{sup}(U)=\mu(e)$ and this belongs to $U$. Also $U$ is clearly a poset, where the order is the one inherited by the usual order in $[0,1]$, and as explained above, $U$ has a greatest element. Each poset $U \subseteq [0, 1]$ which has a greatest element $\alpha$ can be regarded as an inverse monoid with multiplication $\wedge$, that is
$$u \wedge v=\text{min}\{u,v\}.$$
The unit of $(U, \wedge)$ is clearly the element $\alpha$. Now it is obvious that a fuzzy subgroup $(G, \mu , U)$ is a dual premorphism from the group $G$ to the inverse monoid $U$. It is well known that such dual premorphisms are closely related with $F$-inverse covers of inverse monoids. Before we explain what are the $F$-inverse covers of inverse monoids, let us recall a few basic concepts from inverse monoids. An inverse monoid is a monoid $M$ such that for every $x \in M$ there is a unique $x^{-1} \in M$, called the inverse of $x$, such that $xx^{-1}x=x$ and $x^{-1}xx^{-1}=x^{-1}$. Every inverse monoid $M$ comes equipped with a natural partial order $\leq$ defined by: $x \leq y$ if and only if there is an idempotent $e \in M$ such that $x=ye$. Every inverse monoid $M$ has a smallest group congruence which is denoted by $\sigma$ and is characterized by $(x,y) \in \sigma$ if and only if there is an idempotent $e \in M$ such that $xe=ye$. An inverse monoid satisfying the property that each $\sigma$-class contains a greatest element with respect to the natural partial order is called an $F$-inverse monoid. An $F$-inverse monoid $F$ is called an $F$-inverse cover over the group $F/\sigma$ if there exists a surjective idempotent separating homomorphism $F \rightarrow M$. There are a number of important results concerning $F$-inverse monoids which are related with covers and expansions of inverse monoids. The reader can find useful material in papers \cite{ASZ}, \cite{LMS} and \cite{SZ}. Regarding $F$-inverse covers of inverse monoids, as we mentioned before, they are closely related with dual premorphisms between inverse monoids. A dual premorphism $\psi: M \rightarrow N$ between inverse monoids is a map such that $\psi(x^{-1})=(\psi(x))^{-1}$ and $\psi(xy) \geq \psi(x)\psi(y)$ for all $x, y \in M$. The following is Theorem VII.6.11 of \cite{MP} and gives a relationship between dual premorphisms and $F$-inverse covers.
\begin{theorem} \label{R1}
Let $H$ be a group and $M$ an inverse monoid. If $\psi: H \rightarrow M$ is a dual premorphism such that for every $u \in M$, there is an $h \in H$ with $u \leq \psi(h)$, then
$$F=\{(u,h)\in M \times H| u \leq \psi(h)\},$$
is an $F$-inverse cover of $M$ over $H$. Conversely, every $F$-inverse cover of $M$ over $H$ can be so constructed (up to isomorphism).
\end{theorem}
Finally, for everything unexplained here on inverse semigroups we refer the reader to the monographs \cite{ML} and \cite{MP}. While for basics on fuzzy sets and fuzzy groups we refer the reader to \cite{FG}, \cite{R} and \cite{Z}. The book of Mac Lane \cite{MacLane} contains the necessary material on categories and functors,

\section{The definitions of $\mathfrak{F}\mathfrak{G}$ and $\mathfrak{F}\mathfrak{C}$}

We define the category of fuzzy subgroups $\mathfrak{F}\mathfrak{G}$ in the following way. The objects of $\mathfrak{F}\mathfrak{G}$ are triples $(G, \mu, U)$ as defined above, and if $(G, \mu_{1}, U)$ and $(H, \mu_{2}, V)$ are two such triples, a morphism from $(G, \mu_{1}, U)$ to $(H, \mu_{2}, V)$ is a pair $(f, \lambda)$ where $f: G \rightarrow H$ is a group homomorphism, and $\lambda: U \rightarrow V$ is an order preserving map with the property that $\lambda(\text{sup}(U))=\text{sup}(V)$, and, for all $x \in G$
$$\mu_{2}f(x) = \lambda \mu_{1}(x).$$
The unit morphism on an object $(G, \mu , U)$ is defined to be the pair $(1_{G},1_{U})$. Now if $(K, \mu_{3}, W)$ is another object from $\mathfrak{F}\mathfrak{G}$, and $(g, \lambda'): (H, \mu_{2}, V) \rightarrow (K, \mu_{3}, W)$ is another a morphism, we define the composition $(g, \lambda') \circ (f, \lambda)$ as the pair $(gf, \lambda'\lambda)$. We will show that this pair is indeed a morphism from $(G, \mu_{1}, U)$ to $(K, \mu_{3}, W)$. For every $x \in G$, 
\begin{align*}
\mu_{3}gf(x) &= \lambda' \mu_{2}f(x)\\
&= \lambda' \lambda \mu_{1}(x).
\end{align*}
Also $\lambda' \lambda$ is order preserving, and $\lambda' \lambda (\text{sup}(U))= \lambda'(\text{sup}(V))=\text{sup}(W)$. The properties that $\mathfrak{F}\mathfrak{G}$ should satisfy to be a category are straightforward. It looks like, we used in the definition of morphisms of $\mathfrak{F}\mathfrak{G}$ mixed concepts. On the one hand we have homomorphisms of groups, and on the other hand, order preserving maps between posets. In fact, for any two posets $U$ and $V$ with respective greatest elements $\alpha$ and $\beta$, any order preserving map $\lambda:U \rightarrow V$ which sends $\alpha$ to $\beta$, is in fact a homomorphism between monoids $(U, \wedge)$ and $(V, \wedge)$. Indeed, let $u, v \in U$ such that $u \leq v$, then
\begin{align*}
\lambda( u \wedge v)&= \lambda(u) && \text{(since $u \leq v$)}\\
&= \lambda(u) \wedge \lambda(v)  && \text{(since $\lambda(u) \leq \lambda(v)$).}
\end{align*}
In addition to that, the fact that $\lambda(\alpha)=\beta$ says that $\lambda$ is a homomorphism of monoids. The converse is also true, that is, any homomorphism of monoids $\lambda: U \rightarrow V$, is an order preserving map since for every $u,v \in [ 0 , \alpha]$ such that $u \leq v$,
\begin{align*}
\lambda(u)&=\lambda(u \wedge v)&& \text{(since $u \leq v$)}\\
&= \lambda(u) \wedge \lambda(v) && \text{(since $\lambda$ is a homomorphism),}
\end{align*}
which implies that $\lambda(u) \leq \lambda(v)$. The condition that $\lambda(\alpha)=\beta$ follows from the fact that $\lambda$ is a monoid homomorphism. Finally, we remark that $\lambda$ maps the greatest element $\alpha$ of the single $\sigma$-class of $U$ to the greatest element $\beta$ of the single $\sigma$-class of $V$.	

The definition of the category $\mathfrak{F}\mathfrak{C}$ of $F$-inverse covers of inverse monoids is an extension of the definition of the category $\mathfrak{F}$ of $F$-inverse semigroups made in \cite{SZ}. The objects of $\mathfrak{F}\mathfrak{C}$ are triples $(T, M, \varphi)$ where $T$ is an $F$-inverse monoid, $M$ is an inverse monoid, and $\varphi$ is a homomorphism of monoids which is surjective and idempotent separating. We say that $T$ is an $F$-inverse cover of $M$ over $T/\sigma$. If now $(T', M', \varphi')$ is another triple as above, then a morphism from $(T, M, \varphi)$ to $(T', M', \varphi')$ is a pair $(f_{\ast}, \lambda)$ with $f_{\ast}: T \rightarrow T'$ and $\lambda: M \rightarrow M'$ monoid morphisms which map the greatest element of a $\sigma$-class onto the greatest element of some $\sigma$-class, and that satisfy the commutativity condition $\varphi' f_{\ast}=\lambda \varphi$. The identity morphism on the object $(T, M, \varphi)$ is defined the pair $(1_{T}, 1_{M})$ which clearly satisfies the above commutativity condition. The composition of morphisms is defined in the following fashion. If $(f_{\ast}, \lambda): (T,M,\varphi) \rightarrow (T',M',\varphi')$ and $(f'_{\ast}, \lambda'): (T',M',\varphi') \rightarrow (T'',M'',\varphi'')$ are two morphisms, then their composition is defined to be the pair $(f'_{\ast}f_{\ast}, \lambda'\lambda)$. This is indeed a morphism from $(T,M,\varphi)$ to $(T'',M'',\varphi'')$ since
\begin{align*}
\varphi''(f'_{\ast}f_{\ast})&= (\varphi'' f'_{\ast})f_{\ast}\\
&=(\lambda'\varphi')f_{\ast}=\lambda'(\varphi' f_{\ast})\\
&=\lambda'(\lambda \varphi)=(\lambda' \lambda) \varphi,
\end{align*}
and that both compositions $f'_{\ast}f_{\ast}$ and $\lambda'\lambda$ map the greatest element of a $\sigma$-class onto the greatest element of some $\sigma$-class since their respective components do so. Finally, it is easy to see that $\mathfrak{F}\mathfrak{C}$ is indeed a category.

\section{The embedding}

Looking back to the definition of an object $(G, \mu, U)$ from $\mathfrak{F}\mathfrak{G}$, but with $U$ regarded now as an inverse monoid, we see that the map $\mu: G \rightarrow U$ is nothing but a dual premorphism between inverse monoids satisfying the property that for every $u \in U$, there exists $x \in G$ such that $u \leq \mu(x)$. It follows from Theorem \ref{R1}, that there is an $F$-inverse cover of $U$ over $G$ which we write with the long notation $\mathfrak{C}(G, \mu, U)$. More explicitly, 
$$\mathfrak{C}(G, \mu, U)=\{(u,x) \in U \times G| u \leq \mu(x)\}$$
is an inverse monoid whose idempotents turn out to be all the pairs $(u,1)$, where $1$ is the unit of $G$, in particular the unit element is $(\mu(1),1)$. The natural order has a simple description: $(u,x) \leq (v,y)$ if and only if $y=x$ and $u \leq v$. The $\sigma$-class of an element $(u,x)$ consists of all the elements $(v, x)$ with $v \leq \mu(x)$ and its greatest element with respect to the natural order is $(\mu(x),x)$. Finally note that the projection in the first coordinate $\varphi: \mathfrak{C}(G, \mu, U) \rightarrow U$, $(u,x) \mapsto u$ is surjective and idempotent separating. We call the triple $(\mathfrak{C}(G, \mu, U), \varphi, U)$ the $F$-inverse cover associated with the fuzzy subgroup $(G, \mu, U)$. The monoid $\mathfrak{C}(G, \mu, U)$ seems to be useful in connecting inverse semigroups with fuzzy subgroups. An argument which goes in favor to this is that the $\mathcal{H}$-classes of $\mathfrak{C}(G, \mu, U)$ correspond in a way that will be made precise below, to the so called level subsets of $(G, \mu, U)$. Level subsets are defined in \cite{Das} as follows. Given a fuzzy subgroup $(G, \mu, U)$ and $u \in U$, then the level subset $\mu_{u}$ of the fuzzy subset $\mu$ is defined by 
$$\mu_{u}=\{h \in G| \mu(h) \geq u\}.$$
It is proved in Theorem 2.1 of \cite{Das} that such subsets are in fact subgroups of $G$. Before we see the connection they have with the $\mathcal{H}$-classes of $\mathfrak{C}(G, \mu, U)$, we note that $\mathfrak{C}(G, \mu, U)$ is a Clifford monoid. Indeed, it is inverse and its idempotents are central. To see the latter, let $(u,1)$ be an idempotent, and $(v,h)$ an arbitrary element, then
$$(u,1)(v,h)=(u \wedge v, h)=(v \wedge u, h)=(v,h)(u,1).$$ 
To see what an $\mathcal{H}$-class looks like, we recall first that the relations $\mathcal{H}$ and $\mathcal{R}$ coincide in Clifford semigroups. Let now $(v,h) \in \mathfrak{C}(G, \mu, U)$ be such that $(v,h) \mathcal{R} (u,1)$ where $(u,1)$ is some idempotent. There are $(w,a), (w',b) \in \mathfrak{C}(G, \mu, U)$ such that 
$$(v \wedge w, ha)=(v,h)(w,a)=(u,1)$$
and
$$(u \wedge w', bh)=(u,1)(w', b)=(v,h),$$
which both imply that $u=v$. Therefore, if $(v,h) \in H_{(u,1)}$, then necessarily $v=u$. Conversely, any $(u,h) \in \mathfrak{C}(G, \mu, U)$ is $\mathcal{R}$ (hence $\mathcal{H}$)-equivalent with $(u,1)$. Indeed, 
$$(u,h)(u,h^{-1})=(u,1) \text{ and } (u,1)(u,h)=(u,h).$$
All we said means that for any fixed $u \in U$, $H_{(u,1)}=\{(u, h) | u \leq \mu(h)\}$ and this forms a subgroup of $\mathfrak{C}(G, \mu, U)$. Now we show that each level subgroup $\mu_{u}$ is in fact isomorphic to $H_{(u,1)}$. Indeed, the map 
$$\phi: H_{(u,1)} \rightarrow \mu_{u} \text{ such that } (u,h) \mapsto h,$$
is clearly bijective and a homomorphism.

Now we prove our main result.

\begin{theorem}
There is a full and faithful embedding of the category $\mathfrak{F}\mathfrak{G}$ of fuzzy subgroups into the category $\mathfrak{F}\mathfrak{C}$ of $F$-inverse covers of inverse monoids.
\end{theorem}
\begin{proof}
Define $\Omega: \mathfrak{F}\mathfrak{G} \rightarrow \mathfrak{F}\mathfrak{C}$ on objects by sending each fuzzy subgroup $(G, \mu, U)$ to its corresponding $F$-inverse cover $(\mathfrak{C}(G, \mu, U), \varphi, U)$. Further, for each morphism $(f, \lambda): (G, \mu_{1}, U) \rightarrow (H, \mu_{2}, V)$ in $\mathfrak{F}\mathfrak{G}$, if $(\mathfrak{C}(G, \mu, U), \varphi, U)$ and $(\mathfrak{C}(H, \mu', V), \varphi', V)$ are the corresponding $F$-inverse covers, we define
$$f_{\ast}: \mathfrak{C}(G, \mu, U) \rightarrow \mathfrak{C}(H, \mu', V)$$
by setting
$$f_{\ast}(u,x)= (\lambda(u), f(x)).$$
This map is correct since $\lambda(u) \leq \mu_{2}(f(x))$. Indeed, from the definition of the morphism $(f, \lambda)$, we see that 
\begin{align*}
\mu_{2}f(x) &= \lambda \mu_{1}(x)\\
&\geq \lambda(u) && \text{(since $\mu_{1}(x) \geq u$)}.
\end{align*}
Also $f_{\ast}$ is a monoid homomorphism since if $(u,x), (v,y) \in \mathfrak{C}(G, \mu, U)$ such that $u \leq v$, then
\begin{align*}
f_{\ast}((u,x)(v,y))&=f_{\ast}(u \wedge v, xy)\\
&= f_{\ast}(u, xy)\\
&= (\lambda(u), f(xy))\\
&=(\lambda(u) \wedge \lambda(v), f(x)f(y)) && \text{(since $\lambda(u) \leq \lambda(v)$)}\\
&=(\lambda(u), f(x))(\lambda(v), f(y))\\
&=f_{\ast}(u,x)f_{\ast}(v,y).
\end{align*}
Also if $\alpha, \beta$ are the respective units of $U,V$, and $e_{1},e_{2}$ the units of $G, H$ respectively, then
$$f_{\ast}(\alpha, e)=(\lambda(\alpha),f(e_{1}))=(\beta,e_{2}),$$	
so that $f_{\ast}$ preserves the unit element. Further, letting $\varphi_{1}$ and $\varphi_{2}$ be the corresponding maps of the above covers, we see that for every $(u,x) \in \mathfrak{C}(G, \mu, U),$
\begin{align*}
\varphi_{2}f_{\ast}(u,x)&=\varphi_{2}(\lambda(u),f(x))\\
&=\lambda(u)\\
&=\lambda\varphi_{1}(u,x).
\end{align*}
Lastly, if $(\mu(x),x)$ is the greatest element of its $\sigma$-class, then
$$f_{\ast}(\mu(x),x)=(\lambda\mu(x),f(x))=(\mu'f(x),f(x))$$
where $(\mu'f(x),f(x))$ is the greatest element of its $\sigma$-class. Since in addition to what we said, $\lambda$ is a homomorphism of inverse monoids that maps the greatest element $\alpha$ of the only $\sigma$-class of $U$ to the greatest element $\beta$ of the only $\sigma$-class of $V$, then it follows that the pair $\Omega(f, \lambda)=(f_{\ast}, \lambda)$ is a morphism from $(\mathfrak{C}(G, \mu, U), \varphi, U)$ to $(\mathfrak{C}(H, \mu', V), \varphi', V)$. Next we show that $\Omega$ is functorial. It is obvious that when $(f,\lambda)=(1_{G},1_{U})$ is the identity on $(G, \mu, U)$, then $\Omega(1_{G},1_{U})=id_{(\mathfrak{C}(G, \mu, U), \varphi, U)}$. Let now 
$$(f, \lambda): (G, \mu_{1}, U) \rightarrow (H, \mu_{2}, V)$$ 
and 
$$(f', \lambda'): (H, \mu_{2}, V) \rightarrow (K, \mu_{3}, W)$$ 
be two morphisms in $\mathfrak{F}\mathfrak{G}$, and 
$$(f'f, \lambda' \lambda): (G, \mu_{1}, U) \rightarrow (K, \mu_{3}, W)$$ 
their composition, and want to prove that 
$$\Omega(f'f, \lambda'\lambda)= \Omega(f', \lambda')\Omega(f, \lambda),$$
or equivalently that
$$((f'f)_{\ast}, \lambda'\lambda)=(f'_{\ast}, \lambda')(f_{\ast}, \lambda).$$
This is the same as to prove that $(f'f)_{\ast}=f'_{\ast}f_{\ast}$. The latter is true since for every $(u,x) \in \mathfrak{C}(G, \mu_{1}, U)$ we have that
\begin{align*}
(f'f)_{\ast}(u,x)&= ((\lambda' \lambda)(u),(f'f)(x) )\\
&=(\lambda'(\lambda(u)), f'(f(x)))\\
&=f'_{\ast}(\lambda(u),f(x))\\
&=f'_{\ast}(f(u,x)).
\end{align*}
Next we prove that $\Omega$ is faithful. Let $(G, \mu_{1}, U)$ and $(H, \mu_{2}, V)$ be two objects in $\mathfrak{F}\mathfrak{G}$ and 
$$(f, \lambda), (f', \lambda'): (G, \mu_{1}, U) \rightarrow (H, \mu_{2}, V)$$ 
be two parallel morphisms, and assume that $\Omega(f, \lambda)=\Omega(f', \lambda')$. Then, from the definition of $\Omega$, $(f_{\ast}, \lambda)=(f'_{\ast}, \lambda')$, consequently $\lambda=\lambda'$ and $f_{\ast}=f'_{\ast}$. The second equality implies that for every $(u,x) \in \mathfrak{C}(G, \mu_{1}, U)$ we have that $f_{\ast}(u,x)=f'_{\ast}(u,x)$.
It follows that $(\lambda(u), f(x))=(\lambda'(u), f'(x))$, consequently $f(x)=f'(x)$. 

Finally we prove that $\Omega$ is full. Let again $(G, \mu_{1}, U)$ and $(H, \mu_{2}, V)$ be two objects in $\mathfrak{F}\mathfrak{G}$ and 
$$(g, \lambda): (\mathfrak{C}(G, \mu_{1}, U), \varphi_{1}, U) \rightarrow (\mathfrak{C}(H, \mu_{2}, V), \varphi_{2}, V)$$
be a morphism from $\Omega(G, \mu_{1}, U)$ to $\Omega(H, \mu_{2}, V)$. We show that $g$ induces a homomorphism $f: G \rightarrow H$ such that $g=f_{\ast}$ and $\mu_{2}f=\lambda\mu_{1}$. This would prove that $(f,\lambda): (G, \mu_{1}, U) \rightarrow (H, \mu_{2}, V)$ is a morphism in $\mathfrak{F}\mathfrak{G}$ such that $(g,\lambda)=\Omega(f,\lambda)$. For every $x \in G$, let $(\mu_{1}(x),x)$ be the greatest element of its $\sigma$-class in $\mathfrak{C}(G, \mu_{1}, U)$, and let $(\mu_{2}(x'),x')=g(\mu_{1}(x),x)$ which is, from the assumption on $g$, the greatest element of its $\sigma$-class in $\mathfrak{C}(H, \mu_{2}, V)$. It follows that
\begin{align} \label{mu}
\mu_{2}(x')=\varphi_{2}(\mu_{2}(x'),x')=\varphi_{2}g(\mu_{1}(x),x)=\lambda \varphi_{1}(\mu_{1}(x),x)=\lambda\mu_{1}(x).
\end{align}
Since $g$ preserves $\sigma$-classes, there is an induced homomorphism 
$$\tilde{f}: \mathfrak{C}(G, \mu_{1}, U)/\sigma  \rightarrow  \mathfrak{C}(H, \mu_{2}, V)/\sigma$$
which maps the $\sigma$-class $[(\mu_{1}(x),x)]$ to the $\sigma$-class $[(\mu_{2}(x'),x')]$. Considering now the ismorphisms
$$\gamma: G \rightarrow \mathfrak{C}(G, \mu_{1}, U)/\sigma \text{ such that } x \mapsto [(\mu_{1}(x),x)],$$
and
$$\kappa: H \rightarrow \mathfrak{C}(H, \mu_{2}, V)/\sigma \text{ such that } y \mapsto [(\mu_{2}(y),y)],$$
we obtain a homomorphism 
$$f=\kappa^{-1}\tilde{f}\gamma:G \rightarrow H$$
such that $f(x)=x'$ where $x'$ is determined as above. Using (\ref{mu}) we see that $\mu_{2}f(x)=\lambda\mu_{1}(x)$, hence $(f,\lambda)$ is a morphism in $\mathfrak{F}\mathfrak{G}$ from $(G, \mu_{1}, U)$ to $(H, \mu_{2}, V)$. Now we prove that $(g, \lambda)=\Omega(f, \lambda)=(f_{\ast}, \lambda)$ which amounts to saying that $g=f_{\ast}$. Before we prove this, we observe that for every $(u,x) \in \mathfrak{C}(G, \mu_{1}, U)$, the second coordinate of $g(u,x)$ is $x'=f(x)$ as determined above, since $g$ preserves $\sigma$-classes, while the first coordinate is
$$\varphi_{2}g(u,x)=\lambda\varphi_{1}(u,x)=\lambda(u).$$
So 
$$g(u,x)=(\lambda(u), f(x))=f_{\ast}(u, x),$$
consequently, $g=f_{\ast}$ as claimed. This completes the proof.
\end{proof}


\begin{thebibliography}{99}
	

\bibitem{ASZ} K. Auinger, M.B. Szendrei, \textit{On $F$-inverse covers of inverse monoids}, J. Pure Appl. Algebra, 204 (2006), 493-506

\bibitem{ML} M.V. Lawson, \textit{Inverse Semigroup. The Theory of Partial Symmetries}, World Scientific, 1998

\bibitem{LMS} M.V. Lawson, S. W. Margolis, B. Steinberg, \textit{Expansions of inverse semigroups}, J. Aust. Math. Soc. 80, (2006), 205-228

\bibitem{MacLane}  S. Mac Lane, \textit{Categories for the working mathematician}, Springer, 1997

\bibitem{FG} J.N. Mordeson, K.R. Bhutani, A. Rosenfeld, \textit{Fuzzy Group Theory}, Springer, 2005
	

\bibitem{MP} M. Petrich, \textit{Inverse Semigroups}, Wiley, New York, 1984


\bibitem{R} A. Rosenfeld, \textit{Fuzzy groups}, J. Math. Anal. Appl. 35 (1971), 512-517

\bibitem{Das} P. Sivaramakrishna Das, \textit{Fuzzy Groups and Level Subgroups}, J. Math. Anal. Appl. 84 (1981), 264-269 

\bibitem{SZ} M.B. Szendrei, \textit{A note on Birget-Rhodes expansions of groups}, J. Pure Appl. Algebra, 58 (1989), 93-99


\bibitem{Z} L.A. Zadeh, \textit{Fuzzy sets}, Inform. and Comput. 8 (1965), 338-353


	
	
\end{thebibliography}
\end{document}